\documentclass[12pt,a4paper]{article}

\usepackage{a4,enumerate}
\usepackage{amsthm}
\usepackage{latexsym}

\newdimen\templaenge

\DeclareMathAlphabet{\doba}{U}{msb}{m}{n}
\gdef\mC{\doba{C}}

\gdef\mN{\doba{N}}

\gdef\mR{\doba{R}}

\gdef\mZ{\doba{Z}}

\def\qed{{\leavevmode\unskip\nobreak\hfil\penalty 50\hskip 1em%
  \hbox{}\nobreak\hfil\lower 1pt\hbox{$\Box$\kern-.5pt}\parfillskip 0pt
  \finalhyphendemerits 0\par\bigbreak}}
\def\qedmath#1{\setbox0\hbox{$\displaystyle #1$}\templaenge=\textwidth\advance\templaenge by -\wd0%
\setbox1\hbox{$\Box$}\advance\templaenge by -2\wd1%
$$#1\hbox to0pt{\kern.5\templaenge$\Box$\kern-.5pt\hss}$$\par\bigbreak}

\def\al{{\alpha}}
\def\be{{\beta}}

\def\om{{\omega}}

\def\la{{\lambda}}
\def\ka{{\kappa}}

\def\Si{{\Sigma}}
\def\ga{{\gamma}}
\def\ep{{\varepsilon}}
\def\Ga{{\Gamma}}

\def\th{{\vartheta}}
\def\Th{{\Theta}}

\def\phi{{\varphi}}

\def\na{{\nabla}}

\def\el{{\ell}}

\def\cL{\mathcal{L}}

\def\ohne{-}
\def\ti{\tilde}
\def\witi{\widetilde}

\def\nelem{p}
\def\ndim{b}
\def\famindex{n}

\def\halbplatz{\thinspace}
\def\znp{{\setbox0\hbox{\ }\kern.7\wd0}}
\def\vnp{{\setbox0\hbox{\ }\kern\wd0}}

\def\ie{i.\halbplatz e.\ \ignorespaces}

\def\nummerarray#1#2{\par\noindent\setbox0\hbox{\rm (#1)}\setbox1\hbox{$#2$}\unhcopy0%
\dimen0=.5\textwidth \advance\dimen0 by -\wd0 \advance\dimen0 by -.5\wd1 \kern\dimen0 \unhcopy1}

\def\CP{\mC{\mathrm P}}

\def\killvek{K}
\def\psp{P_{\mbox{\rm \scriptsize Spin}}}
\def\pso{{P_{\mbox{\rm \scriptsize SO}}}}

\def\SO{\mathop{{\rm SO}}}

\def\Spin{\mathop{{\rm Spin}}}

\def\Cl{\mathop{{\mC \rm l}}}

\def\tilnz{\hbox{$\witi L_0$}^{\lower1mm\hbox{$\scriptstyle 2$}}}
\def\ker{\mathop{{\rm ker}}}

\def\grad{{\mathop{{\rm grad}}}}

\def\dim{\mathop{{\rm dim}}}

\def\min{\mathop{{\rm min}}}

\def\End{{\mathop{{\rm End}}}}

\def\res#1#2{{#1}\lower .11ex\hbox{$|$}\lower .644ex\hbox{$\scriptstyle #2$}}
\def\stelle#1#2{\left. {#1}\right|_{#2}}

\def\Dhor{{D_{\rm h}^\famindex}}
\def\Dvert{{D_{\rm v}^\famindex}}

\long\def\komment#1{}
\def\beweis#1{{\par\medbreak\noindent {\bf Beweis\setbox0\hbox{#1}%
\ifdim\wd0=0pt .\else\ \ignorespaces #1.\fi}\enspace}}

\def\proof#1{{\par\medbreak\noindent {\bf Proof\setbox0\hbox{#1}%
\ifdim\wd0=0pt .\else\ \ignorespaces #1.\fi}\enspace}}
\def\iop#1{{\par\medbreak\noindent {\bf Idea of proof\setbox0\hbox{#1}%
\ifdim\wd0=0pt .\else\ \ignorespaces #1.\fi}\enspace}}
\def\examples{{\noindent {\bf Examples. }\par\kern-\baselineskip}}

\newtheoremstyle{bemerkungen}{3pt}{3pt}{}{}{\bfseries}{}{ }{}

\newtheorem{theorem}{\bf T{\footnotesize HEOREM}}[section]

\newtheorem{lemma}[theorem]{\bf L{\footnotesize EMMA}}

\theoremstyle{definition}
\newtheorem*{remark}{Remark}

\newtheorem*{example}{Example}

\theoremstyle{bemerkungen}

\begin{document}
\title{The Dirac operator on collapsing $S^1$-bundles}
\author{Bernd Ammann}
\date{October 1998}
\maketitle

\begin{abstract}
We study the behavior of the spectrum of the Dirac operator 
on collapsing $S^1$-bundles. Convergent eigenvalues will 
exist if and only if the spin structure is projectable.


{\bf Mathematics Classification:}
58G25, 58G30, 53C25
\end{abstract}

\section{Introduction}

In this paper we study the spectrum of the Dirac operator 
on collapsing $S^1$-bundles. 

There are some nice results about the behavior of the Laplace operator 
acting on functions 
on a family of collapsing manifolds.
Fukaya \cite{fuk} proved that if a family $(M_\famindex,\ti g_\famindex)_{\famindex\in \mN}$ 
of Riemannian manifolds with 
sectional curvature bounded form above and below and bounded diameter 
converges in the measured Gromov-Hausdorff topology
to a Riemannian manifold $(N,g)$ of lower dimension, 
then the eigenvalues of the Laplace operator acting on functions on 
$(M_\famindex,\ti g_\famindex)$ converge to those of $(N,g)$. 

Fukaya conjectured that it should be possible to replace the bound on the sectional curvature 
by a lower bound on the Ricci curvature. So it seems that the connection between
Gromov-Hausdorff topology and the spectrum should be much closer than we know until today.  

Other papers connecting metric and topological properties to the 
behavior of the spectrum of the 
Laplace operator acting on functions 
are \cite{coco}, \cite{codo} and \cite{wu}.

For the Laplace operator acting on $p$-forms with $p\geq 1$ 
the situation is more complicated.
Until now for $p$-forms there is no analogue to Fukaya's result. But if 
$\pi:(M,\ti g)\to (N,g)$ is a Riemannian submersion whose fibers are minimal submanifolds and 
if the horizontal distribution is integrable, then according to Gilkey and Park \cite{gipa1,gipa2}
the spectrum of the $p$-form Laplacian on $N$\/ is contained in the 
spectrum of the $p$-form Laplacian on $M$. 
The pullbacks of eigenforms on $N$ are eigenforms on $M$ 
to the same eigenvalues.
So if we have a family of such submersions with fixed $(N,g)$,
then this result trivially implies the convergence of certain eigenvalues. 

In a recent paper \cite{glp} Gilkey, Leahy and Park generalized
some weaker version of this result to arbitrary Riemannian submersions.
If an eigenform on $N$ to the eigenvalue $\la$ pulls
back to an eigenform on $M$ to the eigenvalue $\mu$, then $\la\leq \mu$.

In this paper we will concentrate on another elliptic differential operator, the Dirac operator.
We will describe the behavior of the spectrum of the Dirac operator when the fibers of a principal 
$S^1$-bundle collapse. While the Laplacian acting on functions and on forms only depends on the Riemannian metric,
the Dirac operator also depends on the spin structure. It will turn out that for some
spin structures, the projectable ones, 
there are convergent eigenvalues whereas for other spin structures all eigenvalues diverge. 

The methods of this article generalize the methods of \cite{amba} where we only treated the case of
geodesic fibers. There are even families that do not have bounded curvature but for which our results imply
convergence or non-convergence.
Our method is based on a splitting of the Dirac operator into a horizontal Dirac operator, a vertical Dirac 
operator and a zero order term. An analogous 
splitting of the Laplace operator 
has been used in \cite{bebebou}. 

Unfortunately the generalization to non-geodesic fibers has a drawback: 
unlike in the geodesic case 
\cite{amba} we do not get any information about the signs of the
non-convergent eigenvalues of the Dirac operator.

In the proofs of this paper I omitted some technical calculations and I only dealt 
with the case when the dimension $\ndim$ 
of the base manifold is even. More details and the modifications for the odd-dimensional 
case can be found in my thesis \cite{ammanndiss}. I~thank my supervisor Christian B\"ar for 
many interesting discussions and good ideas how to present the collapse result.  
I also want to thank Bruno Colbois who invited me to Chamb\'ery and Grenoble in November 1997.
He told me much about the behavior of the eigenvalues of the Laplacian 
when the manifolds collapse and much about small 
eigenvalues of the Laplace operator.

\section{Spin structures on $S^1$-bundles}

We assume that $S^1$ acts freely and isometrically on 
a compact, connected, oriented $(\ndim +1)$-dimensional Riemannian manifold $(M,\ti g)$, $\ndim \geq 1$. 
Then $M$ is the total space of a principal $S^1$-bundle over some base space $N^\ndim :=M^{\ndim +1}/S^1$.
There is a unique metric $g$ on $N$ such that $\pi :(M,\ti g) \to (N,g)$ is a Riemannian submersion.

The principal $S^1$-bundle $M\to N$ carries a unique connection-1-form 
  $$i\om: TM \to i \mR,$$ 
such that 
$\ker \res{\om}{m}$ is perpendicular to the fibers for any $m\in M$.

The $S^1$-action induces a Killing vector field $\killvek$. 
The fibers of $M\to N$ are (totally) geodesic if and only if $\el:=|\killvek|$ is constant in $m\in M$.
The length of a fiber is $2 \pi \el$.
The metric $\ti g$ on $M$ is completely characterized by $\om$, $\el$ 
and $g$.

Now we recall the notion of a spin structure in order to fix the notation.
The bundle of positively oriented orthonormal frames $\pso(N)$ is a principal $\SO(\ndim)$-bundle.   
The unique non-trivial double covering of $\SO(\ndim)$ will be denoted $\Th: \Spin(\ndim)\to \SO(\ndim)$.
A pair $(\psp(N),\th)$ will be called {\em spin structure} 
if $\psp(N)$ is a principal $\Spin(\ndim)$-bundle over $N$ and if 
$\th:\psp(N)\to\pso(N)$ is a $\Th$-equivariant, fiber preserving map.
Two spin structures $(\psp^1(N),\th_1)$ and $(\psp^2(N),\th_2)$ are isomorphic if there is
a $\Spin(\ndim)$-equivariant, fiber preserving 
isomorphism $A:\psp^1(N)\to\psp^2(N)$ such that $\th_1=\th_2\circ A$.
 
A spin structure exists if and only if the second Stiefel-Whitney class $w_2(TN)$ vanishes. 
In general the spin structure is not unique.
Having chosen a spin structure we can define a spinor bundle on $N$,
Clifford multiplication and a Dirac operator.
The spectrum of the Dirac operator depends on the choice of spin structure.

The definitions of spin structure, Dirac operator, \dots\  on $M$ are clearly analogous.
The $S^1$-action on $M$ induces an $S^1$-action on $\pso(M)$.
A spin structure $\ti\th:\psp(M) \to \pso(M)$ is said to be {\em projectable}
if this $S^1$-action lifts continuously to $\psp(M)$.

\examples
\begin{enumerate}[(1)]
\item We view the $(\ndim +1)$-dimensional 
torus $T^{\ndim +1}$ as the total space of the circle bundle 
$T^{\ndim+1}\to T^{\ndim}$. The torus 
$T^{\ndim +1}$ carries $2^{\ndim +1}$ spin structures, half of them are projectable, 
half of them are non-projectable. 
Recently the special case $T^2\to T^1=S^1$ with the trivial 
spin structure and non-geodesic fibers has been intensively studied
in \cite{fried2}.
\item We view $S^{2l+1}$ as the total space of the Hopf fibration $S^{2l+1}\to \mC P^l$.
The (unique) spin structure on $S^{2l+1}$ is projectable 
if and only if $\mC P^l$ is spin and therefore
if and only if $l$ is odd.
\end{enumerate}

There is a natural isomorphism from projectable spin structures on $M$ to spin structures on $N$.
However $M$ may admit a non-projectable spin structure, even if there are no spin structures on $N$.
The Hopf fibration $S^{2l+1}\to \CP^{l}$ with $l$ even is an example for this phenomenon.

In this paper we will not only look at one single $S^1$-fibered space $M$, but at a family $M_\famindex $ 
of such spaces over a base space $N$ that does not depend on $\famindex$. 
For each $\famindex \in \mN$ the manifold $M_\famindex $ 
carries a Riemannian metric~$\ti g_\famindex $, such that $\pi_\famindex :(M_\famindex ,\ti g_\famindex )\to (N,g)$ is an 
$S^1$-bundle and a Riemannian submersion. The quantities $\el_\famindex $ and $\om_\famindex $ are defined as they are defined for $M$.
Roughly speaking, we will analyze the behavior of the spectrum of the Dirac operator if $\el_\famindex $ tends to $0$, whereas 
$\om_\famindex $ stays small in a suitable sense. 
This situation will be called {\em collapse}.

Note that $N$ and $g$ do not depend on $\famindex $. On the other hand $M_\famindex $ may even change its topological type for
different values of $\famindex $.

\section{Projectable spin structures}

Let $(N,g)$ be a Riemannian manifold carrying a spin structure that shall be fixed throughout this section.
The eigenvalues of the Dirac operator $D$ on  
$N$ will be denoted by $(\mu_j)_{j\in\mN}$.
Moreover let $(M_\famindex ,\ti g_\famindex )_{\famindex \in \mN}$ be a family of Riemannian manifolds
such that $S^1$ acts freely and isometrically on each $M_\famindex $,
and let $(M_\famindex ,\ti g_\famindex )\to (N=M_\famindex /S^1,g)$ be Riemannian submersions.
 
The spin structure on $M_\famindex $ shall be the unique 
projectable spin structure corresponding to the spin structure on $N$. 
As above let $2\pi \el_\famindex $ be the length of the fibers and $i\om_\famindex $
the connection-1-form on $M_\famindex \to N$.
Furthermore we assume the collapsing condition
\begin{equation}\label{kollapsbed}
\begin{array}{c}
  \|\el_\famindex \cdot d\om_\famindex \|_\infty \to 0  \mbox{ and } \|\el_\famindex \|_\infty \to 0 \mbox{ for } \famindex  \to \infty \\[2mm]
  \al:=\limsup\limits_{\famindex \to \infty}\|\grad \el_\famindex \|_\infty<1 
\end{array}
\end{equation}

\begin{theorem}\label{kollaps}
The eigenvalues $\big(\la_{j,k}(\famindex )\big)_{j\in\mN,k\in\mZ}$ 
of the Dirac operator $\witi D^{\famindex}$ on  
$M_\famindex $ 
can be numbered in such a way that:
\begin{enumerate}[\rm (1)]
\item For all $\ep>0$ there is a $\famindex _0\in \mN$, 
      such that we have for any $\famindex  \geq \famindex _0$ and $j\in\mN$, $k \in \mZ\ohne\{0\}$  
      \begin{eqnarray*}
        \|\el_\famindex \|_\infty^2\;{\la_{j,k}(\famindex )}^2& \geq & |k|\,(|k|-\al) \, - \,\ep.
      \end{eqnarray*}
      In particular we get 
      ${\la_{j,k}(\famindex )}^2\to \infty$ for $\famindex \to \infty$.
      \newline\medskip
      Furthermore, if $M_\famindex $ and $\om_\famindex $ do not depend on $\famindex $, then we also have 
      for $j\in\mN$, $k \in \mZ\ohne\{0\}$.       
      \begin{eqnarray*}
        \limsup_{\famindex \to\infty}  \Big(\min_{\nelem\in N}\el_\famindex (\nelem)\Big)^2\,{\la_{j,k}(\famindex )}^2
        & \leq & |k|\,(|k|+\al).
      \end{eqnarray*}
      This upper bound of ${\la_{j,k}(\famindex )}^2$ is not uniform in $j$ and $k$.
\item If $\ndim=\dim N$ is even, then we get for $\famindex \to \infty$
      $$\la_{j,0}(\famindex ) \to \mu_j.$$
      However for $\ndim=\dim N$ odd we obtain
      \begin{eqnarray*}
         \la_{2j-1,0}(\famindex ) & \to & \phantom{-} \mu_j\\
         \la_{2j,0}(\famindex )   & \to &           - \mu_j
      \end{eqnarray*}
      In both cases the convergence of the eigenvalues $\la_{j,0}(\famindex )$ is uniform in $j$.
\end{enumerate} 
\end{theorem}

\examples
\begin{enumerate}[(1)]
\item The square of the Dirac operator on the flat torus $\mR/a\mZ\times \mR/\el \mZ$ with
the trivial spin structure has the eigenvalues
  $$\la_{j_1,j_2,k}^2 = 4\pi^2 \left( {j_1^2\over a^2} + {k^2\over \el^2}\right),\qquad
    j_1\in \mZ,\; j_2\in\{0,1\},\;k\in\mZ.$$ 
\item We take the standard metric on $\mC P^l$, the $\om$ coming from the Hopf fibration
$S^{2l+1}\to \mC P^l$
and a sequence of constant functions $(\el_n)_{n\in\mN}$, $\el_n>0$, $\el_n\to 0$.
The metrics on $S^{2l+1}$ characterized by these quantities are called {\em Berger metrics} 
$g_{\el_n}$.
If $l$ is odd, some eigenvalues diverge whereas other eigenvalues converge 
to the eigenvalues of $\mC P^l$.
The spectrum of the spheres with Berger metrics has been explicitely calculated in \cite{hi,b2}.
C.~B\"ar used the collapse results 
to calculate the spectrum of the Dirac operator on $\mC P^l$ \cite{amba}.
\end{enumerate}

\iop{}
In order to prove Theorem~\ref{kollaps} we write the Dirac operator 
$\witi D^\famindex  $ as a sum of a ``vertical Dirac operator'', a ``horizontal Dirac operator''
and a zero order term. 

For the definition of the vertical Dirac operator we need the Lie derivative of spinors 
in the direction of the Killing field.
The $S^1$-action on $\psp{M_\famindex }$ induces an $S^1$-action 
on $\Si M_\famindex =\psp(M_\famindex )\times_{\Spin(\ndim+1)}\Si_{\ndim+1}$.
The latter action will be denoted by $\ka$.
A spinor with base point $m\in M$ will be mapped by $\ka(e^{is})$ 
to a spinor with base point $m\cdot e^{is}$.
We define the {\em Lie derivative} of a smooth spinor field $\Psi$ 
in the direction of the Killing field $\killvek$ as
  $$\cL_\killvek(\Psi)(m)= \stelle{d \over ds}{s=0}
    \ka(e^{-is})(\Psi(m\cdot e^{is})).$$   
As $\cL_\killvek$ is the differential of a representation of the Lie group $S^1$ 
on $L^2(\Si M_\famindex )$, we get the decomposition
  $$L^2(\Si M_\famindex )=\bigoplus_{k \in \mZ} V_{k,\famindex }$$
into the eigenspaces $V_{k,\famindex }$ of the operator $\cL_\killvek$ to the eigenvalue $ik$, $k\in \mZ$.
The $S^1$-action commutes with the Dirac operator on $M_\famindex $,
and therefore this decomposition is respected by the Dirac operator.

We calculate the difference between the covariant derivative 
and the Lie derivative 
in the direction of $\killvek$.
For any smooth section $\ti\Psi$ of $\Si M_\famindex $ we get 
\begin{equation}\label{vertnabla}
\na_\killvek \ti\Psi -  \cL_\killvek \ti\Psi= 
    {\el_\famindex ^2 \over 4}\, \ga(d\om_\famindex ) \ti\Psi 
    - {1\over 2}\ga(\killvek/\el_\famindex )\ga(\grad \el_\famindex )\ti\Psi.
\end{equation}
Here $\ga(V)$ resp.\  $\ga(\be)$ denotes Clifford multiplication of a 
spinor by the vector $V$ and the $2$-form $\be$ resp.

What we got until now is some kind of Fourier decomposition along the fibers.
In the following it will turn out that for any $k$ there is a natural 
isomorphism $Q_{k,\famindex }$ from the summand $V_{k,\famindex }$ 
to the vector space of sections of a twisted spinor bundle on the quotient space $N$.

To the  principal $S^1$-bundle $\pi_\famindex :M_\famindex\to N$ we associate the complex vector bundle
$L_\famindex :=M_\famindex \times_{S^1}\mC$ with a connection given by $i \om_\famindex $. 
The horizontal lift of a vector (field) $X$ on $N$ to a vector (field) on 
$M_\famindex$ will be denoted by 
$\witi X$.  

For $\ndim$ even, the isomorphism $Q_{k,\famindex }$ is given by the following technical 
lemma that is proven in
\cite[Lemma-Definition 7.2.3]{ammanndiss}. 

\begin{lemma}\label{horinabla}
Let $\ndim=\dim N$ be even.
Then for any $\famindex $ there is an isometry of Hilbert spaces 
  $$Q_{k,\famindex }:L^2(\Si N \otimes L_\famindex ^{-k}) \to  V_{k,\famindex }$$
such that the horizontal covariant derivative is given by
  $$\na_{\witi X}{Q_{k,\famindex }(\Psi)} =Q_{k,\famindex }(\na_{X}\Psi) 
    + {\el_\famindex  \over 4}\, \ga(\killvek/\el_\famindex  )\ga(\witi W_X){Q_{k,\famindex }(\Psi)} -
      {X(\el_\famindex )\over 2\el_\famindex } Q_{k,\famindex }(\Psi),$$
where $W_X$  is the vector field on $N$ satisfying 
$d\om_\famindex (\witi X,\cdot)=\langle \witi W_X, \cdot \rangle$.
Clifford multipliction is preserved, \ie    
   $$Q_{k,\famindex }(\ga(X)\Psi)=\ga(\witi X) {Q_{k,\famindex }(\Psi)}.$$
\end{lemma} 

Now we are ready to prove the theorem for $\ndim$ even.
We define the {\it horizontal Dirac operator}
as the unique closed linear operator 
$\Dhor :L^2(\Si M_\famindex ) \to L^2(\Si M_\famindex )$ on each $V_{k,\famindex }$ given by 
  $$\Dhor := Q_{k,\famindex }\circ D^\famindex  \circ {Q_{k,\famindex }}^{-1}$$
where $D^\famindex $ is the twisted Dirac operator on
$\Si N \otimes L_\famindex ^{-k}$.

We define the {\it vertical Dirac operator} 
  $$\Dvert:=\ga(\killvek/\el_\famindex )\,\cL_\killvek$$ 
and the zero order term
  $$Z_\famindex :=-(1/4)\,\ga(\killvek/\el_\famindex  )\,\ga(d\om_\famindex ).$$
 
Using Formula (\ref{vertnabla}) and Lemma~\ref{horinabla}
we can express the Dirac operator as a sum:
\begin{eqnarray*}
   \witi D^\famindex  & = & {1\over \el_\famindex  } \Dvert + \Dhor  + \el_\famindex   Z_\famindex .\\
\end{eqnarray*}

Since $\Dhor $, $\ga(\killvek/\el_\famindex )$ and $Z_\famindex $ commute with the $S^1$-action
they also commute with $\cL_\killvek$.
Therefore each summand of the above decomposition of the Dirac operator 
maps each $V_{k,\famindex }$ into itself.
Furthermore the spectrum of $\res{\Dhor}{V_{0,\famindex }}$ is just the
spectrum of the Dirac operator $D$ acting on sections of $\Si N$. 
So the eigenvalues of $\res{\witi D^\famindex }{V_{0,\famindex }}$
converge to the eigenvalues of $D$. 
This fact immediately implies (2) of the theorem 
for $\ndim$ even.

An elementary calculation shows that $\Dhor$ anticommutes with 
$\ga(\killvek/\el_\famindex )$ and therefore it anticommutes with
$\Dvert:=\ga(\killvek/\el_\famindex )\cL_\killvek$.

From this anticommutativity we get the formula
  $$\Dhor \left({1\over \el_\famindex } \Dvert\right)+\left({1\over \el_\famindex }\Dvert\right)\Dhor =
     \ga\Big(\grad (1/\el_\famindex )\Big)\Dvert.$$
We take the square of $A^\famindex :=(1/\el_\famindex )\Dvert+\Dhor$ and get
  $$(A^\famindex )^2= {1\over \el_\famindex ^2}\left(\Dvert\right)^2+ \left(\Dhor\right) ^2 - 
    {\ga(\grad \,\el_\famindex )\over \el_\famindex ^2}\Dvert.$$
Now for $\al_\famindex :=\|\grad\, \el_\famindex \|_\infty$ and for  $\ti\Psi\in V_{k,\famindex }$, $k\neq 0$ we have
\begin{eqnarray*}
\left( (A^\famindex )^2 \ti\Psi,\ti\Psi\right)_{M_\famindex } & 
= & \left(\left({k^2\over \el_\famindex ^2} + \left(\Dhor\right)^2 
-ik {\ga(\grad \,\el_\famindex )\ga(\killvek/\el_\famindex )\over \el_\famindex ^2}\right)\ti\Psi,\ti\Psi\right)_{M_\famindex }\\
& \geq & {|k|(|k|-\al_\famindex )\over \|\el_\famindex \|_\infty^2}\left(\ti\Psi,\ti\Psi\right)_{M_\famindex }.
\end{eqnarray*} 

The eigenvalues of $\res{(A^\famindex )^2}{V_{k,\famindex }}$ are therefore larger than or 
equal to
  $$\Bigl(|k|\,(|k|-\al_\famindex )\Bigr)/\|\el_\famindex \|_\infty^2.$$ 
As $\el_\famindex \cdot d\om_\famindex $ is a bounded operator, whose norm vanishes for $\famindex \to \infty$,
the norm of the zero order term tends to zero. So we get the first part of (1) of the theorem.

Now let $M_\famindex $ and $\om_\famindex $ be independent from $\famindex $. Then also $L_\famindex $ and $D^\famindex $ are independent from~$\famindex $.
Analogously we obtain
\begin{eqnarray*}
\left( (A^\famindex )^2 Q_{k,\famindex }(\Psi),Q_{k,\famindex }(\Psi)\right)_{M_\famindex } & 
\leq  & \left( {(D^\famindex )}^2\,\Psi,\Psi\right)_{N}\\[1mm]
&&{}+ {|k|(|k|+\al_\famindex)\over \Big(\min_{\nelem\in N} \el_\famindex (\nelem)\Big)^2}\left(Q_{k,\famindex }(\Psi),Q_{k,\famindex }(\Psi)\right)_{M_\famindex }\\
\end{eqnarray*}
From this we get the second part of (1) for $\ndim$ even.

The proof for $\ndim$ odd runs quite analogously, but we have to do some modifications 
as in this case the rank of the spinor bundle over $M_\famindex $ is two times the rank of the spinor bundle 
over $N$. For details about this case we refer to \cite{amba} or \cite{ammanndiss}.
\qed

\begin{remark}
Suppose we have a complex vector bundle $E\to N$ with a metric connection $\na^E$. 
Now we replace the Dirac operator $D$ by the twisted Dirac operator $D^E$ acting 
on sections of $\Si N \otimes E \to N$, the Dirac operator $\witi D^\famindex $ will be replaced by 
the twisted Dirac operator $\witi D^{\famindex ,E}$ acting on sections of 
$\Si M_\famindex  \otimes \pi_\famindex ^* E \to M_\famindex $. Then Theorems~\ref{kollaps}, \ref{kollapsnp} and \ref{abschnp}
are still valid. The proof of this generalized version is essentially the same as the proof above. 
\end{remark}

\section{Non-projectable spin structures}

Now we turn to the case of non-projectable spin structures. 
In this case we get a similar result, but with a stronger restriction on the gradient of the fiber length.
The variable $k$ from the last section does
no longer take integer values, but values in $\mZ+(1/2)$. 
Therefore all eigenvalues will diverge.

We define $N$, $g$, $M_\famindex $, $\ti g_\famindex $, $\om_\famindex $ and $\el_\famindex $ as above. 
However in contrast to the last section we do not assume that $N$ carries any spin structure.
Instead we suppose that each $M_\famindex $ carries a non-projectable 
spin structure.

The collapsing condition in this section is
\begin{equation}\label{kollapsbednp}
\begin{array}{c}
  \|\el_\famindex \cdot d\om_\famindex \|_\infty \to 0  \mbox{ and } \|\el_\famindex \|_\infty \to 0 \mbox{ for } \famindex  \to \infty \\[2mm]
  \al:=\limsup_{\famindex \to \infty}\|\grad\, \el_\famindex \|_\infty<1/2\,.
\end{array}
\end{equation}

\begin{theorem}\label{kollapsnp} 
All eigenvalues of the Dirac operator $\witi D^\famindex $
on $M_\famindex $ diverge. 

Moreover, the eigenvalues $(\la_{j,k}(\famindex ))_{j\in\mN,k\in\mZ+(1/2)}$ can be numbered in such a way that 
\begin{enumerate}[\rm (a)]
\item For any $\ep>0$ there is $\famindex _0\in \mN$, 
      such that for any $\famindex  \geq \famindex _0$ and $j\in\mN$, $k \in \mZ+(1/2)$ we get
      \begin{eqnarray*}
        \|\el_\famindex \|_\infty^2\;{\la_{j,k}(\famindex )}^2& \geq & |k|\,(|k|-\al) \, - \,\ep.
      \end{eqnarray*}
\item If $M_\famindex $ and $\om_\famindex $ are independent from $\famindex $, then we additionally have for 
      any $j\in\mN$, $k \in \mZ+(1/2)$       
      \begin{eqnarray*}
        \limsup_{\famindex \to\infty}  \left[\Big(\min_{\nelem\in N}\el_\famindex (\nelem)\Big)^2\,{\la_{j,k}(\famindex )}^2\right]
        & \leq & |k|\,(|k|+\al).
      \end{eqnarray*}
      This upper bound is not uniform in $j$ and $k$.
\end{enumerate}
\end{theorem}

\begin{example}
Suppose we have a family of collapsing Berger metrics on $S^{2l+1}$ 
as in the previous example,
but with $l$ even. Then all eigenvalues diverge.
The same holds for tori with non-projectable spin structures.
\end{example}

\iop{}
The proof is a variation of the proof of Theorem \ref{kollaps}.
We will restrict to the case $\ndim$ even.

Let $\ti\th_\famindex :\psp(M_\famindex )\to \pso(M_\famindex )$ be a non-projectable spin structure on $M_\famindex $. 
In this case $N$ may or may not be spin. 
We define $P_{\SO(\ndim)}(M_\famindex )$ to be the set of all frames in
$\pso(M_\famindex )$ having $K/\el_\famindex $ as first vector. Then $P_{\SO(\ndim)}(M_\famindex )$
is a principal $\SO(\ndim)$-bundle over $M_\famindex$. 
Moreover $P:=\ti\th_\famindex ^{-1}(P_{\SO(\ndim)}(M_\famindex ))$ is a principal $\Spin(\ndim)$-bundle 
over $M_\famindex $. 

The action of $S^1\cong (\mR / 2\pi \mZ)$ does not lift to $P$, 
but the double covering of $S^1$, i.e. $S^1\cong (\mR/ 4\pi \mZ)$,
does act on $P$.
We define $\Spin^{\mC}(\ndim)$ to be $\Spin(\ndim)\times_{\mZ_2}(\mR/4\pi \mZ)$ 
where $-1\in\mZ_2$ identifies $(-A,c)$ with $(A,c+2\pi)$.
The complex standard representation $\Si_\ndim$ of $\Spin (\ndim)$ is also
a representation for $\Spin^{\mC}(\ndim)$ where $s+4\pi \mZ\in \mR/4\pi\mZ$
operates as $\exp( i s/2)$.  
The actions of $\Spin(\ndim)$ and $\mR/4\pi\mZ$ on $P$ induce a free action of 
$\Spin^{\mC}(\ndim)$ on $P$ and we can view $P$ as a principal 
$\Spin^{\mC}(\ndim)$-bundle over $N$. Then we can form the bundle 
$P\times_{\Spin^{\mC}(\ndim)} \Si_\ndim$. 

If $N$ is spin, this bundle is just $\Si N \otimes L_\famindex ^{1\over 2}$. 
If $N$ is not spin, then neither $\Si N$ nor $L_\famindex ^{1\over 2}$ exist
but $\Si N \otimes L_\famindex ^{1\over 2}$ does exist.

Again we get a splitting 
  $$L^2(\Si M_\famindex )=\bigoplus_{k\in \mZ+{1\over 2}} V_{k,\famindex }$$
into eigenspaces $V_{k,\famindex }$ for $\cL_K$ to the eigenvalue $ik$. 
The rest of the proof of Theorem \ref{kollapsnp} is the same as
the one for the case $k\neq 0$ in Theorem \ref{kollaps}.
\qed

The formulas in the proof also show the following lower bound for eigenvalues of $D^2$.
We define the {\em Clifford norm} $\|\eta\|_{\Cl}$ of a 2-form $\eta$ on $M$ to be the $L^2$-operator norm 
of $\ga(\eta)\in\End(\Ga\Si M)$.

\begin{theorem}\label{abschnp}  
Let $(M,\ti g)$ be a Riemannian manifold on which  
$S^1$ acts freely and isometrically. 
The quotient $N:=M/S^1$ shall carry the unique metric $g$ for which 
$M\to N$ is a Riemannian submersion. 
We assume that $M$ carries a non-projectable spin structure.

If $\al:=\|\grad\,\el\|_\infty<1/2$, then all eigenvalues $(\la_i)_{i\in\mZ}$
of the Dirac operator $\witi D$ on $M$
satisfy
  $$|\la_i|\geq {\sqrt{1-2\,\|\grad\, \el\|_\infty}\over 2\|\el\|_\infty}
     -{\|\el\,d\om\|_{\Cl}\over 4}.$$        
\end{theorem} 

\begin{remark}
This lower bound for the smallest eigenvalue of the Dirac operator depends on the spin structure.

The existence of such a bound is not surprising: on many standard spaces
the spectrum of the Dirac operator is explicitely 
known, among many others 
flat tori \cite{fried1} and compact quotients
of $3$-dimensional Heisenberg groups with a left invariant metric \cite{amba}. 
For these examples it turns out that the spectrum
strongly depends on the spin structure, in particular the smallest eigenvalue
also depends on it. So for metrics close enough to these standard spaces
there should be a lower bound for the
smallest eigenvalue of the Dirac operator depending on the spin structure.

Until recently such a bound was not known. Theorem~\ref{abschnp} 
provides such a bound for the case of a circle bundle with non-projectable
spin structure and sufficiently short fibers.
If $\|\grad\,\el\|_\infty$ and $\|\el\|_\infty$ are sufficiently small Theorem~\ref{abschnp}
gives a better bound than any other estimate published before.
For more details and other recent estimates depending on the spin structure I refer 
to my thesis \cite{ammanndiss}.
\end{remark}

\end{document}